\documentclass[10pt]{amsart}

\usepackage{amsmath,amssymb,amsfonts,epsfig}

\newcommand{\rb}{\raisebox}
\newcommand{\ig}{\includegraphics}
\newcommand\risS[6]{\rb{#1pt}[#5pt][#6pt]{\begin{picture}(#4,15)(0,0)
  \put(0,0){\ig[width=#4pt]{#2.eps}} #3
     \end{picture}}}
\newcommand{\chd}[1]{\risS{-9}{#1}{}{25}{20}{15}}
\newtheorem{thm}{Theorem}[section]
\newtheorem{defn}[thm]{Definition}
\newtheorem{Example}[thm]{Example}
\newcommand{\e}{\varepsilon}
\newcommand{\ph}{\varphi}
\newcommand{\Hom}{\mathop{\rm Hom}\nolimits}
\newcommand{\sign}{\mathop{\rm sign}}
\newcommand{\im}{\mathop{\rm im}}
\def\cro{{\mathsf x}}
\def\ar{{\mathsf a}}
\def\kr{{\mathfrak{C}}}
\def\ard#1{\risS{-12}{#1}{}{25}{15}{17}}
\renewcommand{\leq}{\leqslant}
\renewcommand{\geq}{\geqslant}

\begin{document}


\title[Polyak-Viro formulas for the Conway polynomial]{Polyak-Viro formulas for 
coefficients of the Conway polynomial}

\author{SERGEI CHMUTOV, MICHAEL CAP KHOURY, and ALFRED ROSSI}

\address{The Ohio State University, Mansfield,
1680 University Drive, Mansfield, OH 44906.
{\tt chmutov@math.ohio-state.edu}\linebreak
Department of Mathematics, The University of Michigan,
530 Church St. Ann Arbor, \mbox{MI 48109.
{\tt mjkhoury@umich.edu}} \linebreak
Department of Physics, The Ohio State University, 
191 W Woodruff Ave., Columbus, OH 43210.
{\tt rossi.49@osu.edu}}

\maketitle

\begin{abstract}
We describe the Polyak-Viro arrow diagram formulas for the coefficients of the Conway polynomial. As a consequence, we obtain the Conway 
polynomial as a state sum over some subsets of the crossings of the knot diagram. It turns out to be a simplification of a special case of Jaeger's state model for the HOMFLY polynomial.
\end{abstract}



\section*{Introduction} \label{s:intro}

In this paper we are working with the Conway polynomial $\nabla(L)$ of an oriented link $L$ defined by the equations
\newcommand{\unkn}{\rb{-4.2mm}{\ig[width=10mm]{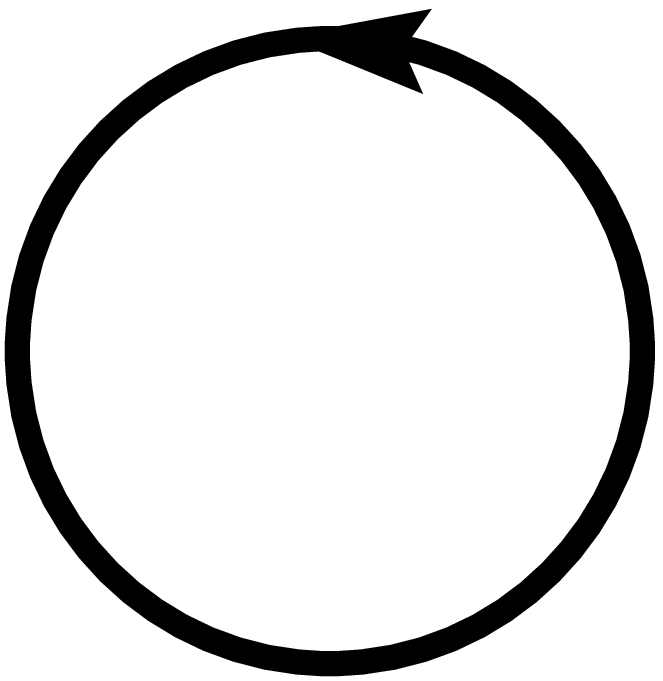}}}
\newcommand{\rlints}{\rb{-4.2mm}{\ig[width=10mm]{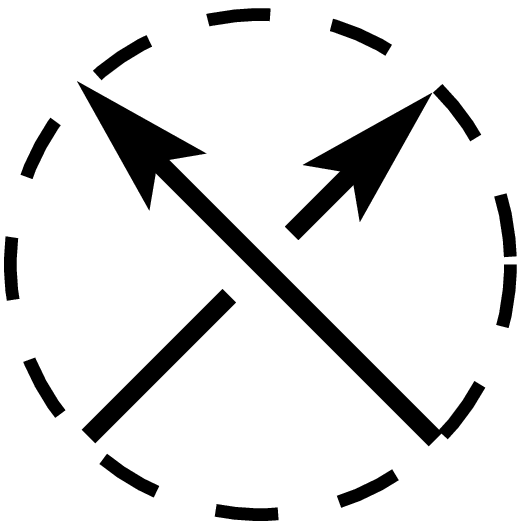}}}
\newcommand{\lrints}{\rb{-4.2mm}{\ig[width=10mm]{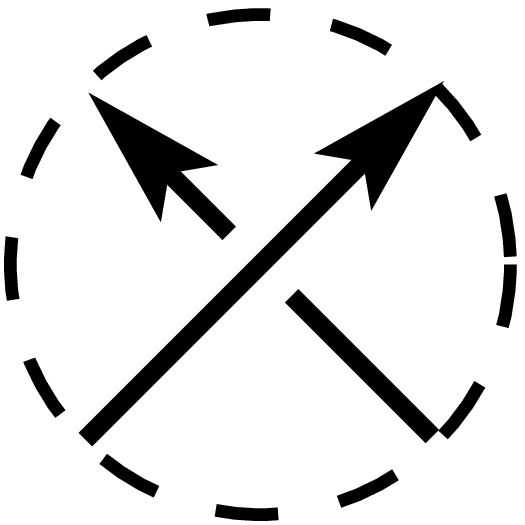}}}
\newcommand{\twoup}{\rb{-4.2mm}{\ig[width=10mm]{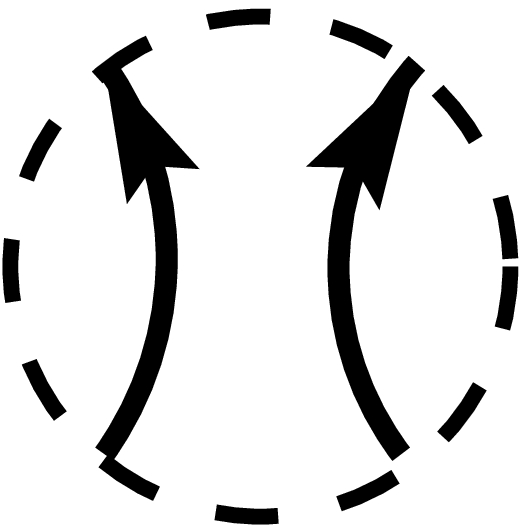}}}
$$
\nabla\Bigl(\lrints\Bigr)-\nabla\Bigl(\rlints\Bigr) = z  
    \nabla\Bigl(\twoup\Bigr)\ , \hspace{1cm}
\nabla\Bigl(\unkn\Bigr) = 1\ .
$$
Its coefficient $c_n(L)$ at $z^n$ is a Vassiliev invariant of order
$\leq n$. The purpose of this paper is to provide, in case $L$ is a knot, some formulas for $c_n(L)$ in terms of certain subdiagrams of the Gauss diagram of $L$. These formulas may be equivalently reformulated (Section \ref{s:st-model}) as a state model on a diagram of the knot. L.~Kauffman noted that this state model should be related to Jaeger's state model for the HOMFLY polynomial \cite{Ja}. Indeed it turns out that our state model is a simplification of Jaeger's model to the special case of knots and to the Conway polynomial. Also our formulas lead to two (different) extensions of the Conway polynomial to long virtual knots.

In Section \ref{s:st-model} we formulate the state model for the Conway
polynomial. We review Gauss diagrams and Polyak-Viro formulas in Section
\ref{s:gd-pv}. In Section \ref{s:mth} we formulate our main result (Theorem
\ref{th: main}) in terms of Gauss diagrams and give some examples. Section
\ref{s:proof} is devoted to the proof of the main theorem.

\section{State model} \label{s:st-model}

A subset $S$ of the crossings of a knot diagram $K$ is said to be
{\it one-component} if the curve obtained from $K$ by smoothing all the  crossings of $S$ according to orientation has one component.

Assume that the diagram $K$ has a base point and $S$ is a one-component subset of the crossings. 
Let us travel along the smoothened curve starting at the base point. In this journey we pass a neighborhood of every smoothing twice.
We call the subset $S$ {\it ascending} if, for every smoothing, the first time we approach its neighborhood on overpass of $K$ (so we jump down to perform the smoothing) and upon returning to the neighborhood we approach it on the underpass (jumping up).

Define the {\it down} polynomial, in variable $z$, as
$$\nabla_{\mbox{\scriptsize asc}}(K)\quad := 
\sum_{\genfrac{}{}{0pt}{2}{S \mbox{\scriptsize\quad ascending}
     }{\mbox{\scriptsize one-component}}}
   \Bigl( \prod_{{\mathsf x}\,\in S} \mbox{wr}(\cro) \Bigr)\ z^{|S|}\ ,
$$
where $\mbox{wr}(\cro)$ is the local writhe of the crossing $\cro$. If $S$ is the empty set, then we set the product to be equal to 1 by definition. Therefore the free
term of $\nabla_{\mbox{\scriptsize asc}}(K)$ always equals 1.

\def\plinc{\mbox{\begin{picture}(5,5)(0,0)
          \put(2.5,2){\circle{8}} \put(0,0){\mbox{\scriptsize $+$}}
                 \end{picture}}}
\def\miinc{\mbox{\begin{picture}(5,5)(0,0)
          \put(2.5,2){\circle{8}} \put(-.5,0){\scriptsize \mbox{$-$}}
                 \end{picture}}}
For example, for the knot $6_2$ 
$$\risS{-32}{k6-2or}{\put(-50,30){\mbox{\tt knot $6_2$}}}{70}{35}{40}
\qquad\qquad
G_{6_2}=\ \risS{-25}{Gd6-2}{\put(65,30){\mbox{\tt Gauss diagram}}
        \put(90,20){\mbox{\tt of $6_2$}}
             }{60}{20}{25}$$
there are eleven one-component subsets with two crossings,
$\{1,2\}$, $\{1,4\}$, $\{1,5\}$, $\{1,6\}$, $\{2,4\}$, $\{2,5\}$, $\{2,6\}$, 
$\{3,4\}$, $\{3,5\}$, $\{4,6\}$, $\{5,6\}$. However, only three of these subsets, $\{2,4\}$, $\{2,6\}$, and $\{4,6\}$, are ascending. The products of the writhes for the subsets $\{2,5\}$, $\{2,6\}$, and $\{4,6\}$ are equal to $-1$, $+1$, and $-1$, respectively. Hence the coefficient of $z^2$ in the polynomial 
$\nabla_{\mbox{\scriptsize asc}}(6_2)$ equals $-1+1-1=-1$.

\smallskip
{\bf Corollary of Theorem \ref{th: main}.} {\it The Conway polynomial $\nabla(K)$ of a knot $K$ is equal to the down polynomial of its diagram,
$$\nabla(K) = \nabla_{\mbox{\scriptsize asc}}(K)\ .$$}

Let us remind that the Conway polynomial of the knot $6_2$ is equal to 
$\nabla(6_2)=1-z^2-z^4$. So indeed its coefficient at $z^2$ equals $-1$.

\smallskip
This formula also holds for links. 
Similarly to $\nabla_{\mbox{\scriptsize asc}}(K)$, one can define the {\it descending} polynomial $\nabla_{\mbox{\scriptsize des}}(K)$. It turns out that for all classical knots
$\nabla_{\mbox{\scriptsize asc}}(K) = \nabla_{\mbox{\scriptsize des}}(K)$. This equality fails for virtual knots. However, each of these polynomials can be extended to virtual links with a based point. See more details in Remark \ref{fin-rem}.2.

\section{Gauss diagrams and Polyak-Viro formulas}\label{s:gd-pv}

\begin{defn}\rm
A {\it Gauss diagram} is a chord diagram with oriented chords and with
numbers $+1$ or $-1$ assigned to each chord. 
\end{defn}

With a knot diagram we associate a Gauss diagram with the outer circle being the parameterizing circle $S^1$ of our knot, a chord for each double point of the diagram, each chord oriented from the
overpass to the underpass and the local writhe number assigned to each
double point (chord). An example of the Gauss diagram of the knot $6_2$ is given in the introduction.

The Reidemeister moves for knot diagrams can be expressed in terms of Gauss diagrams
as follows (see, for example, \cite{CDbook}).
$$\begin{array}{c@{\qquad}c}
\mbox{R-I}_{\mbox{Gd}}:&
\risS{-20}{virrI}{
            \put(24,15){\mbox{$\scriptstyle \e$}}
            \put(186,14){\mbox{$\scriptstyle \e$}}
             }{200}{20}{35} \\
\mbox{R-II}_{\mbox{Gd}}:&
\risS{-20}{virrII1}{
            \put(10,18){\mbox{$\scriptstyle \e$}}
            \put(24,18){\mbox{$\scriptstyle -\e$}}
             }{120}{20}{35}\\
\mbox{R-III}_{\mbox{Gd}}:&
\risS{-18}{virrIII4}{}{120}{20}{25}\qquad\quad
\risS{-18}{virrIII5}{}{120}{0}{0}\ .
\end{array}$$

M.~Polyak and O.~Viro suggested \cite{PV} the following approach to 
represent knot invariants in terms of Gauss diagrams.

\begin{defn}\rm
An {\it arrow diagram} is a based chord diagram with oriented chords. 
\end{defn}

\begin{defn}\rm
Let $A$ be an arrow diagram and let $G$ be a Gauss diagram, both with base points.
A {\it homomorphism} $\ph$ from $A$ to $G$, $\ph\in\Hom(A,G)$, is an orientation preserving homeomorphism of the circle of $A$ to the circle of $G$ which maps the base point to the base point and induces an injective map of chords of $A$ to chords of $G$ respecting the orientation of the chords.
\end{defn}

\begin{defn}\rm
The {\it pairing} between a based arrow diagram and a based Gauss diagram 
is defined by
$$\langle A,G\rangle :=
 \sum_{\ph\in\Hom(A,G)}\quad\prod_{c\mbox{\scriptsize\ chord in\ }A}\sign(\ph(c))\ .
$$
\end{defn}

We want to use this pairing to define knot invariants by choosing
an arrow diagram $A$ and then sending 
$K\mapsto G(K)\mapsto \langle A,G(K)\rangle$.
This invariant will be well defined if the result is
independent of the choice of the Gauss diagram $G(K)$ and the base point on it.
For example, this is the case for the arrow diagram 
$A=\risS{-6}{cd22arw}{}{15}{10}{8}$.
If $G$ is the Gauss diagram of the knot $6_2$ from the introduction, then
there are three homomorphisms of $A$ into $G$, which send the two arrows 
of $A$ to the pairs
of chords $\{2,5\}$, $\{2,6\}$, and $\{4,6\}$ of $G$, respectively. Thus, in this case
$\langle A,G\rangle = -1$.

\medskip
In general, if you take an arbitrary arrow diagram $A$,
the value $\langle A,G(K)\rangle$ is not
uniquely defined by the knot $K$.
However, if we extend the pairing to a linear combination of arrow diagrams
$$\langle \sum_i \lambda_i A_i,G \rangle :=
\sum_i \lambda_i \langle A_i,G \rangle
$$
by linearity, then there are many linear combinations of arrow 
diagrams that yield knot invariants by this construction.
Moreover, with a slight generalization of arrow diagrams involving signed arrows, there is a
general theorem due to M.~Goussarov \cite{G,GPV} stating
that any Vassiliev invariant can be obtained from a suitable
linear combination of arrow diagrams (possibly with signed chords).

In the next section we explicitly describe the linear combinations of arrow diagrams which give the coefficients of the Conway polynomial.

The Polyak-Viro formulas may be considered for links as well. For example, the following formula gives the linking number for a link $L$ with two components $K_1$ and $K_2$:
$$lk(K_1,K_2) = 
\langle\,\risS{-10}{adln}{}{65}{15}{20}\,,G(L)\rangle\ .
$$
Such formulas for 2-component links will be used in the proof of our Main theorem.

\section{Main Theorem}\label{s:mth}

\begin{defn}\label{def:one-comp}\rm
A chord diagram $D$ is said to be $k${\it -component} if after parallel doubling of each chord according to the picture 
$\risS{-6}{chord}{}{30}{12}{8}\ 
 \risS{-2}{totor}{}{25}{0}{0}\ \risS{-6}{dchord}{}{30}{0}{0}\ ,
$
the resulting curve will have $k$ components. We use the notation $|D|=k$.
\end{defn}

\begin{Example}\rm 
For chord diagram with two chords we have:
$$\Bigl|\chd{cd22ch4}\Bigr| =1\ \Longleftarrow\ \chd{cd22-ppar}\ ,\hspace{2cm}
  \Bigl|\chd{cd21ch4}\Bigr| =3\ \Longleftarrow\ \chd{symst4}\ .
$$
In this paper we will work with one-component diagrams only.
With four chords, there are four one-component diagrams (the notation is borrowed from \cite{CDbook}): 
$$d^4_1=\chd{cd4-01}\ ,\quad d^4_5=\chd{cd4-05}\ ,\quad 
  d^4_6=\chd{cd4-06}\ ,\quad\mbox{and}\quad d^4_7=\chd{cd4-07}\ .
$$
\end{Example}

\begin{defn}\rm
We can turn a one-component chord diagram with a base point into an arrow diagram according to the following rule.
{\it Starting from the base point we travel along the diagram with doubled chords. During this journey we pass both copies of each chord in opposite directions. Choose an arrow on each chord which corresponds
to the direction of the first passage of the copies of the chord.} Here is an example.\vspace{-5pt}
$$\chd{cd22ch4}\quad  \risS{2}{totor}{}{25}{15}{15}\quad
  \chd{cd22par}\quad  \risS{2}{totor}{}{25}{0}{0}\quad \chd{cd22arw}\ .
$$
\end{defn}
We call the arrow diagram obtained in this way the {\it ascending} arrow diagram.

\begin{defn}\label{def:cc}\rm
The {\it Conway combination} $\kr_{2n}$ is the sum of all based one-component ascending arrow diagrams with $2n$ arrows. For example,
$$\begin{array}{rcl}
\kr_2 &:=& \risS{-12}{cd22arw}{}{25}{15}{20}\ ,\\
\kr_4 &:=& \ard{cd4-01arw}\ 
  +\ \ard{cd4-07arw1} + \ard{cd4-07arw2} + \ard{cd4-07arw3} + \ard{cd4-07arw4}+ \\   
&&\hspace{-8pt} 
   + \ard{cd4-05arw1} + \ard{cd4-05arw2} + \ard{cd4-05arw3} + \ard{cd4-05arw4} 
   + \ard{cd4-05arw5} + \ard{cd4-05arw6} + \ard{cd4-05arw7} + \ard{cd4-05arw8} + \\   
&&\hspace{-8pt} 
   + \ard{cd4-06arw1} + \ard{cd4-06arw2} + \ard{cd4-06arw3} + \ard{cd4-06arw4} 
   + \ard{cd4-06arw5} + \ard{cd4-06arw6} + \ard{cd4-06arw7} + \ard{cd4-06arw8}\ .
\end{array}$$
Note that for a given one-component chord diagram we have to consider all possible choices for the base point. However, some choices may lead to the same arrow diagram. In
$\kr_{2n}$ we list them without repetitions. For instance, all choices of a base point for the diagram $d^4_1$ give the same arrow diagram. So $d^4_1$ contributes only one arrow diagram to $\kr_4$. The diagram $d^4_7$ contributes four arrow diagrams because of its symmetry, while $d^4_5$ and $d^4_6$ contribute eight arrow diagrams each.
\end{defn}

\begin{thm}\label{th: main}
For $n\geqslant 1$, the coefficient $c_{2n}$ of $z^{2n}$ in the Conway polynomial 
of a knot $K$ with the Gauss diagram $G$ is equal to
$$c_{2n} = \langle \kr_{2n},G \rangle\ .$$ 
\end{thm}

\begin{Example}\label{ex:g-2-eva}\rm
Consider the knot $K:=6_2$ and its Gauss diagram $G:=G_{6_2}$ from the introduction. To compute the pairing $\langle \kr_4,G\rangle$ we have to match
the arrows of each diagram of $\kr_4$ with the arrows of $G$. One common property of the arrows in $\kr_{2n}$ is that the first (and the last) arrow end-point we meet while traveling along the circle counterclockwise (starting with the base point) is the tail of the arrow. This follows from the above arrow rule for $\kr_{2n}$. Hence the arrow $\{1\}$ of $G$ can not participate in the matching with any diagram of $\kr_4$. The only candidates to match with the first arrow of a diagram of $\kr_4$ are the arrows $\{2\}$ and $\{4\}$ of $G$. If it would be $\{4\}$, then 
$\{1, 2, 3\}$ do not participate in the matching, and there would remain only 3 arrows to match with the four arrows of $\kr_4$. Therefore the arrow of $G$ which matches with the first arrow of a diagram of $\kr_4$ must be $\{2\}$. In a similar way we can find that the arrow of $G$ which matches with the last arrow of a diagram of $\kr_4$ must be $\{6\}$. This leaves three possibilities to match with the 
four arrows of $\kr_4$: $\{2,3,4,6\}$, $\{2,3,5,6\}$, and $\{2,4,5,6\}$. Checking them all we find only one quadruple, $\{2,3,5,6\}$, which matches with the second diagram of the second row of $\kr_4$. The product of the local writhes of the arrows $\{2,3,5,6\}$ is equal to
$(-1)(-1)(+1)(-1)=-1$. In other words,
$$\langle \kr_4,G \rangle = 
  \langle\ \ard{cd4-05arw2}\ ,G \rangle = -1\ ,$$
which coincides with the coefficient $c_4$ of the Conway polynomial 
$\nabla(K)= 1-z^2-z^4$.
\end{Example}

\section{Proof}\label{s:proof}

Let us regard $\langle \kr_{2n},G_D \rangle$ as a function of the knot diagram $D$.
The proof consists of two parts. In the first part we study how the function
changes with a switching of a crossing of $D$ (exchanging the over-strand and the under-strand
at the crossing) and produce a skein relation for our invariant which models the 
Conway skein relation. This would involve two-component links and an extension of 
the function to their diagrams. In the second part we prove the coincidence with 
the Conway polynomial using an induction on the number of arrows of the Gauss 
diagram $G_D$.

\subsection{Skein relation.}\label{ss:sk-rel}
The notions of a one-component ascending arrow diagram and the 
Conway combination from Definition \ref{def:cc} can be naturally extended 
to (arrow) diagrams with two circles. In this case the number of arrows must 
be odd, so we have:
$$\begin{array}{rcl}
\kr_1 &:=& \risS{-11}{adln}{}{65}{20}{20}\ ,
\end{array}$$
$$\begin{array}{rcl}
\kr_3 &:=& \risS{-11}{c3cr1}{}{65}{20}{20} +
  \risS{-11}{c3br1}{}{65}{0}{0} + \risS{-11}{c3br4}{}{65}{0}{0}
       + \risS{-11}{c3br5}{}{65}{0}{0} \\   
&&\hspace{-8pt} 
 + \risS{-11}{c3ar1}{}{65}{10}{20} + \risS{-11}{c3ar2}{}{65}{0}{0}+
  \risS{-11}{c3ar3}{}{65}{0}{0} + \risS{-11}{c3ar4}{}{65}{0}{0}  \\   
&&\hspace{-8pt} 
 + \risS{-11}{c3ar5}{}{65}{10}{20} + \risS{-11}{c3ar6}{}{65}{0}{0}\ .
\end{array}$$

{\bf Lemma.} {\it Suppose $K_+$ is a knot diagram with a positive distinguished
crossing $\cro$, and $K_-$ and $K_0$ are the corresponding knot and 2-component
link obtained by changing the crossing $\cro$ as in the Conway skein relation:
$$K_+=\risS{-9}{lrints}{\put(10,3){$\cro$}}{25}{15}{10}\hspace{2cm}
  K_-=\risS{-9}{rlints}{\put(10,3){$\cro$}}{25}{0}{0}\hspace{2cm}
  K_0=\risS{-9}{twoup}{}{25}{0}{0}\ .
$$
Let us introduce shorter notations:
$G_+:=G_{K_+}$, $G_-:=G_{K_-}$, and $G_0:=G_{K_0}$.
Then, 
\begin{equation}\label{eq:sk-re-kn}
\langle\kr_{2n},G_+\rangle - \langle\kr_{2n},G_-\rangle
 = \langle\kr_{2n-1},G_0\rangle.
\end{equation}
}

\medskip
{\bf Proof.} 
Let $A$ be one of the arrow diagrams of $\kr_{2n}$. If $\ph_+\in\Hom(A,G_+)$ is 
a homomorphism such that $\cro\not\in\im\ph_+$ then such a homomorphism exists 
in $\Hom(A,G_-)$ as well, and they cancel each other on the left side.
Now suppose for some arrow $\ar\in A$, $\ph_+(\ar)=\cro$. We can construct a two-circle one-component ascending arrow diagram $A_\ar$ and a homomorphism
$\ph_\ar\in\Hom(A_\ar,G_0)$ whose contribution to  
$\langle\kr_{2n-1},G_0\rangle$ is the same:
\begin{equation}\label{eq:hom-contr}
\prod_{c\mbox{\scriptsize\ chord in\ }A}\sign(\ph_+(c))=
\prod_{c\mbox{\scriptsize\ chord in\ }A_\ar}\sign(\ph_\ar(c))\ .
\end{equation}
The arrow diagram $A_\ar$ is obtained from $A$ by doubling the chord
$\ar$ as in the definition \ref{def:one-comp}. It has two circles, and obviously it is one-component and ascending. Also the diagram $A_\ar$ contains $2n-1$ arrows. 
Note that the Gauss diagram $G_0$ of the link $K_0$ is obtained from $G_+$ by a similar doubling of the arrow $\cro$ (more precisely, of the arrow corresponding to the crossing $\cro$). So any homomorphism $\ph_+:A\to G_+$ which sends $\ar$ to $\cro$ induces a homomorphism 
$\ph_\ar:A_\ar\to G_0$ which sends the arrows of $A_\ar$ (which may be identified with the arrows of $A$ different from $\ar$) to the same arrows of $G_0$ (which may be identified with the corresponding arrows of $G_+$). Then the equation (\ref{eq:hom-contr}) is obvious since 
$\sign(\ph_+(\ar))=\sign(\cro)=1$.

In a similar way, a homomorphism $\ph_-:A\to G_-$ which sends some arrow $\ar$to $\cro$ induces a homomorphism $\ph_\ar:A_\ar\to G_0$. However, since $\sign(\ph_-(\ar))=\sign(\cro)=-1$, now the left and right hand sides of
(\ref{eq:hom-contr}) differ by a sign. But the homomorphism
$\ph_-\in\Hom(A,G_-)$, as a part of $\langle\kr_{2n},G_-\rangle$, occurs in
the left hand side of the desired equation of the lemma with the sign $-1$.
Therefore its contribution to the left hand side will be the same as the  contribution
of the corresponding $\ph_\ar$ to the right hand side,
$\langle\kr_{2n-1},G_0\rangle$.

In the other direction for any homomorphism $\ph_0:A_0\to G_0$, where $A_0$ is a two-circle arrow
diagram from $\kr_{2n-1}$, we can construct either $\ph_+$ or $\ph_-$ which
contributes into the left hand side of the equation of the lemma the same amount as $\ph_0$ to $\langle\kr_{2n-1},G_0\rangle$. Indeed, $\ph_0$
maps some two arcs from different circles of $A_0$ to the two arcs of $G_0$ corresponding to the two pieces of $K_0$ in the vicinity of $\cro$ on the picture above. We can make a connected sum of the two circles of $A_0$ by connecting these arcs with a band. The result will be a one-circle arrow diagram. We put an extra chord $\ar$ across the band and orient it in the direction which makes the whole diagram ascending. (This direction depends on which of the two arcs of $A_0$ was passed first while the traveling along $A_0$ with doubled chords). The band corresponds to a
half-twisted band making a connected sum of the components of $K_0$ to produce either $K_+$ or $K_-$, depending on the orientation of the new arrow $\ar$. The resulting $K_+$ or $K_-$ determines the sign of $\ar$. Thus we obtain a one-component ascending arrow diagram $A$ with a
distinguished arrow $\ar$ and its homomorphism, either $\ph_+$ or $\ph_-$. It is easy to see that this construction is indeed inverse to the construction above: $\ph_0=\ph_\ar$. This proves the lemma. \hspace{\fill}$\square$

\medskip 
{\bf Example.} Let us continue our example with the knot $K=6_2$ and its
Gauss diagram $G=G_{6_2}$. 
Let us choose crossing $\{5\}$ as the
distinguished one. Then we can denote the knot diagram $K$ as $K_+$. The corresponding knot $K_-$ is not interesting because there are no homomorphisms from $\kr_4$ to $G_-$ (its Gauss diagram is obtained from $G$ by reversing its arrow $\{5\}$ and changing its sign to
$-1$). The link
$K_0$ is more interesting: 
$$K_0=\ \risS{-32}{k6-2-0}{}{70}{35}{40}\qquad
G_0=\ \risS{-25}{Gd6-2-0}{}{60}{20}{25}\ =\
  \risS{-18}{Gd6-2-0s}{}{85}{0}{0}\ .
$$
Doing an analysis similar to Example \ref{ex:g-2-eva} one can conclude that there in only one non-trivial homomorphism $\ph_0$ from the first arrow diagram $A_0$ of $\kr_3$ to $G_0$ which sends the arrows to $\{2,3,6\}$ of $G_0$. The arcs of $K_0$ in the vicinity of the crossing $\{5\}$ are represented on the Gauss diagram $G_0$ by the parallel copies of the arrow $\{5\}$ from $G$. On the picture of $A_0$, the preimages of these arcs under the homomorphism $\ph_0$ are located in between the two crossing arrows of the left circle and on the right portion of the right circle. These are the places where we are
supposed to make a band connected sum.
$$A_0=\risS{-11}{c3cr1}{}{65}{25}{17}
\quad  \risS{-2}{totor}{}{25}{0}{0}\quad
\risS{-11}{ar1}{}{65}{0}{0}\ =\ 
\risS{-11}{ar2}{\put(30,26){$\ar$}}{65}{0}{0}
\ =\ \risS{-12}{cd4-05arw2}{\put(-5,5){$\ar$}}{25}{0}{0} = A
$$
Also we are supposed to put an arrow $\ar$ across the band. In order to make the diagram $A$ ascending we have to orient it down. Thus we obtain the diagram $A$, the only diagram of $\kr_4$ contributing to $\langle\kr_4,G_{6_2}\rangle$.

\medskip
An analogous lemma may be formulated for two-component links. Namely, {\it let $L_+$, $L_-$, and $L_0$ be a triple of diagrams participating in the Conway skein relation, and $G_+$, $G_-$, $G_0$ be their Gauss diagrams. We assume that $L_+$ and $L_-$ are two-component links, and the two strands at the crossing $\cro$ belong to different components. Then $L_0$ will be a knot diagram. We have
\begin{equation}\label{eq:sk-re-li}
\langle\kr_{2n+1},G_+\rangle - \langle\kr_{2n+1},G_-\rangle
 = \langle\kr_{2n},G_0\rangle.
\end{equation}}

The proof is similar to the proof of the lemma. We use these two skein relations to simplify the diagrams in the next subsection.

\subsection{Coincidence with the Conway coefficients.}
We proceed by induction on the number of arrows of the Gauss diagram $G_D$, where $D$ is either a knot or a 2-component link diagram.

If $D$ has no crossings then there is nothing to prove. 

Now let us assume $\langle \kr_{2n},G_D \rangle=c_{2n}$ and 
$\langle \kr_{2n+1},G_D \rangle=c_{2n+1}$ for all knot (2-component link) diagrams $D$
with less than $m$ crossings. Let $D$ be a diagram with $m$ crossings. 
We can pick a crossing on $D$ and use the relations
(\ref{eq:sk-re-kn}) and (\ref{eq:sk-re-li}) to simplify the corresponding 
Gauss diagrams. The diagrams on the right hand sides of these relations have less than $m$ crossings. So, by induction, the right hand sides coincide with the corresponding Conway coefficients. For dealing with the left hand sides we need to consider the cases of knots and links separately.

{\bf Knots.} Changing the appropriate crossings using the relation 
(\ref{eq:sk-re-kn}) we can make our diagram $D$ {\it descending}. This means that traveling along the knot diagram starting from the base point,  the first passage of each crossing will be an overpass. 
On Gauss diagrams this means that all arrows are oriented in accordance with the orientation of the circle of the Gauss diagram, i.e. traveling along the circle for every arrow we first pass its tail and then its head. 
Hence we can represent $\langle \kr_{2n},G_D \rangle$ as 
$\langle \kr_{2n},G_{D'} \rangle$, for some descending diagram $D'$, plus some terms of the form $\langle \kr_{2n},G_{D_0} \rangle$, where the 2-component link diagram $D_0$ has less than $m$ crossings. The descending diagram $D'$ represents an unknot, so its Conway polynomial is equal to 1. Therefore its coefficients $c_{2n}$ ($n\geq 1$) vanish. On the other hand, any subdiagram of a descending Gauss diagram is also descending. None of the arrow diagrams of $\kr_{2n}$ has this property, so there is no homomorphism of $\kr_{2n}$ to $G_{D'}$ and 
$\langle \kr_{2n},G_{D'} \rangle=0$. By induction 
$\langle \kr_{2n},G_D \rangle=c_{2n}(D)$.

{\bf Links.} Now $D$ represents a  2-component link with $m$ crossings.
Using (\ref{eq:sk-re-li}) we change some crossing between the components of $D$ in order to lower the component with the base point to the bottom. On Gauss diagrams, this means that all arrows between different components will point toward the component with the base point. So we represent
$\langle \kr_{2n+1},G_D \rangle$ as $\langle \kr_{2n+1},G_{D'} \rangle$ for some diagram $D'$ with the base component lower than the other one, plus some terms of the form $\langle \kr_{2n},G_{D_0} \rangle$ where the knot diagram $D_0$ has less than $m$ crossings. For the link $D'$ with unlinked components the Conway polynomial is equal to zero. Also in $\kr_{2n+1}$ there are no diagrams whose arrows are all oriented toward the based component. Thus $\langle \kr_{2n+1},G_{D'} \rangle=0$ and  
$\langle \kr_{2n+1},G_D \rangle=c_{2n+1}(D)$ follows by induction.

This completes the proof of Main Theorem.
\hspace{\fill}$\square$

\subsection{Remarks.}\label{fin-rem}

{\bf 1.} The equation $\nabla(L) = \nabla_{\mbox{\scriptsize asc}}(L)$ for links follows from the fact that a morphism of a {\em connected} arrow diagram to the Gauss diagram of $L$ can be reconstructed from its image.

{\bf 2.} M.~Polyak found a direct way to prove the invariance of 
$\langle \kr_{2n},G_D \rangle$ under the Reidemeister moves. Together with the skein relation of Section \ref{ss:sk-rel} it gives a direct proof of our Main Thereom. This way also shows that the formulas of Main Theorem \ref{th: main} provide an extension of 
$\nabla_{\mbox{\scriptsize asc}}(L)$ (as well as 
$\nabla_{\mbox{\scriptsize des}}(L)$) to virtual links with a based point. In particular, to long virtual knots. 
The
second coefficients of these two etensions are the two second order invariants of long virtual knots from \cite{GPV}.
The Polyak-Viro formulas for Vassiliev invariants coming from the HOMFLY polynomial were found in \cite{CP}. That paper contain also the extensions
of the HOMFLY polynomial to virtual string links.

\section*{Acknowledgments}
The work has been done as part of the Summer 2006 VIGRE working group\\
\verb#http://www.math.ohio-state.edu/~chmutov/wor-gr-su06/wor-gr.htm#\\
``Knots and Graphs" 
at the Ohio State University, funded by NSF grant DMS-0135308. We are grateful to L.~Kauffman, M.~Polyak, and O.~Viro for valuable discussions. S.Ch.~thanks the Max-Plank-Institut f\"ur Mathematik in Bonn (where this paper was finished) for excellent working conditions and warm hospitality.

\end{document}